\documentclass[11pt]{amsart}

\usepackage{amsmath} %Advanced math typing
\usepackage{amsthm} %For theorem styles
\usepackage{amsfonts} %For \mathbb
\usepackage{amssymb} %For more math symbols
\usepackage{stmaryrd} %For even more math symbols (like \llbracket)
\usepackage{mathtools} %For extensible math symbols
\usepackage{thm-restate} %For restatable theorems

\usepackage[style=alphabetic, maxnames=100, maxalphanames=100]{biblatex} %For references
\DeclareFieldFormat{postnote}{\mknormrange{#1}} %Prevents biblatex from automatically adding p.

\usepackage{xifthen} %For \ifthenelse

\usepackage{xcolor} %For colored text
\usepackage{bm} %For bold symbols with \bm
\usepackage{bbm} %For alternative blackboardsymbols with \mathbbm

\usepackage{array} %For special commands in tables
\usepackage[inline]{enumitem} %For enumerate styles

\usepackage{microtype} %Corrects kerning
\usepackage{lmodern} %Loads Latin modern fonts (helps microtype package)

\usepackage{pdflscape} %Allows alternating to landscape mode
\usepackage{afterpage} %For \afterpage macro

\usepackage{tikz} %For tikzfigures
\usetikzlibrary{cd} %For commutative diagrams
\usetikzlibrary{calc} %For coordinate calculations

\usepackage[in]{fullpage} %Gets rid of stupidly large margins leaving 1 inch on every side

\usepackage[hypcap=true]{subcaption} % For captions in subfloats

\usepackage[%
  plainpages=false,% Usage of not only arabic numbers
  pdfpagelabels,% Enable page labels
  colorlinks,% Colored links
  citecolor=black,% Citation color
  linkcolor=black,% Link color
  urlcolor=black,% URL color
  filecolor=black,% File color
  bookmarksopen=false% For starting document with bookmarks tree closed
]{hyperref} %For pdf links
\usepackage[all]{hypcap} %Fixes figures and tables anchors for hyperref
\hypersetup{hypertexnames=false}%Fixes hyperref duplicate destinations

\makeatletter
\newcommand*{\newletterthm@internal}{}%Dummy definition to use \renewcommand later
\newcommand*{\newletterthm}[1]{%
  \def\newletterthm@name{#1}%
  \renewcommand*{\newletterthm@internal}[1][]{%
    \ifthenelse{\isempty{##1}}{%
      \expandafter\expandafter\expandafter\newtheorem%
      \expandafter\expandafter\expandafter{%
        \expandafter\newletterthm@name%
        \expandafter}%
      \expandafter{%
        \newletterthm@text}%
      \expandafter\renewcommand%
      \expandafter*%
      \expandafter{%
        \csname the#1\endcsname}{\Alph{#1}}%
    }{%
      \expandafter\expandafter\expandafter\newtheorem%
      \expandafter\expandafter\expandafter{%
        \expandafter\newletterthm@name%
        \expandafter}%
      \expandafter{%
        \newletterthm@text}[##1]%
    \expandafter\renewcommand%
        \expandafter*%
        \expandafter{%
          \csname the#1\endcsname}{\csname the##1\endcsname.\Alph{#1}}%
    }%
  }%
  \newletterthm@newthm%
}

\newcommand*{\newletterthm@newthm}[2][]{%
  \ifthenelse{\isempty{#1}}{%
    \def\newletterthm@text{#2}%
    \newletterthm@internal%
  }{%
    \expandafter\newtheorem\expandafter{\newletterthm@name}[#1]{#2}%
  }%
}
\makeatother
\newtheoremstyle{thmstyle}%style name
  {\medskipamount}%space before
  {\smallskipamount}%space after
  {\slshape}%font used
  {0pt}%indentation
  {\bfseries}%modifier theorem head
  {.}%punctuation between theorem head and body
  { }%space after punctuation
  {\thmname{#1}\thmnumber{ #2}{\normalfont\thmnote{ (#3)}}}%theorem specifier

\newtheoremstyle{plainstyle}%style name
  {\medskipamount}%space before
  {\smallskipamount}%space after
  {\rmfamily}%font used
  {0pt}%indentation
  {\bfseries}%modifier theorem head
  {.}%punctuation between theorem head and body
  { }%space after punctuation
  {\thmname{#1}\thmnumber{ #2}{\normalfont\thmnote{ (#3)}}}%theorem specifier

\theoremstyle{thmstyle}

\newletterthm{theoremLet}{Theorem}
\newletterthm{problemLet}{Problem}

\theoremstyle{plainstyle}

\makeatletter
\def\refdescformat#1{%
  \phantomsection%
  \let\oldlabel\label%
  \let\label\@gobble%
  \edef\@currentlabel{\MakeLowercase{#1}}% Label format (in \ref).
  \let\label\oldlabel%
  #1.% Item format (in description).
}
\makeatother

\newlist{refdesc}{description}{1}
\setlist[refdesc]{format={\refdescformat}}
\newlist{enumdef}{enumerate}{1}
\setlist[enumdef]{before={\leavevmode}, label={\arabic*.}, ref={\thetheorem.\arabic*}}

\setlist[enumerate]{label={\roman*.}, ref={(\roman*)}} % Provides default enumitem arguments for enumerate

\makeatletter
\newcommand{\biggg}{\bBigg@\thr@@}

\newcommand{\Biggg}{\bBigg@{3.5}}

\newcommand{\bigggg}{\bBigg@{4}}

\newcommand{\Bigggg}{\bBigg@{4.5}}

\newcommand{\biggggg}{\bBigg@{5}}

\newcommand{\Biggggg}{\bBigg@{5.5}}

\makeatother
\numberwithin{equation}{section} % Changes equation numbering to be within section.

\newcommand{\mcf}{\mathcal{F}}

\newcommand{\br}{\vspace{2mm}}

\let\epsilon\varepsilon

\newcommand{\mcp}{\mathcal{P}}
\newcommand{\mch}{\mathcal{H}}

\title{Remarks on a recent preprint of Chernikov and Towsner} 

\author{Maryanthe Malliaris}

\begin{document}

\maketitle

\begin{abstract}
In this brief note, we first give a counterexample to a theorem in Chernikov and Towsner,
2510.02420v1. In 2510.02420v2, the theorem has changed but as we explain the proof has a mistake.  
The change in the statement, due to changes in the underlying definition, affects the paper's claims.  
Since that theorem had been relevant to connecting the work of their 
paper to Coregliano-Malliaris high-arity PAC learning, a connection which 
now disappears, 
we also explain 
why their definitions miss crucial aspects that our work was designed to grapple with. 
\end{abstract}

\br

Earlier this month, Chernikov and Towsner posted a short note [CT25] proposing to 
quickly reprove an equivalence of $k$-dependence with 
a notion of learnability and a notion of packing, 
and so derive some of our results by alternate methods based on their [CT20].

[CT25] has two versions, and since the 
descriptions of the claims have not all been updated to reflect the 
mathematical changes across the two versions, 
it is necessary to compare them to see what was actually done.  
In the first version of their paper, the main theorem 
has an easy counterexample, given below.  
We haven't tracked down where the mistake is, whether in their earlier paper [CT20] 
or in the deduction from that paper.   
In the second version of their paper, a definition in the main theorem 
has changed and so the theorem is different, and of course, has 
different implications, and is no longer related to our work.  
We have not verified whether the revised Chernikov-Towsner theorem is correct,
but we note that their proposed proof contains some elementary mistakes.

\br

\noindent 
For background: 

In the last few years, in a series of papers, Coregliano and Malliaris developed
a working theory called high-arity PAC learning, and gave 
complete characterizations in two different regimes, involving 
appropriate versions of uniform convergence, partite learning, nonpartite learning, agnostic learning,
classical and probabilistic packing lemmas, combinatorial dimension, structured correlation.
One of the two combinatorial dimensions involved was Shelah's $k$-dependence.
The papers were: 
\begin{itemize}
\item[1.] ~[Coregliano-M 24] ``High-arity PAC learning via exchangeability.''
February 2024, 145 pages.
arXiv:2402.14294. (main theorems stated on pps. 51-53)

\item[2.] ~[Coregliano-M 25a] ``A packing lemma for $VCN_k$-dimension and learning
high-dimensional data.” May 2025, 29 pages. arXiv:2505.15688. 

\item[3.] ~[Coregliano-M 25b] ``Sample completion, structured correlation and Netflix problems.''
September 2025, 97 pages.  arXiv:2509.20404. (main theorems stated on pps. 38-40)
\end{itemize}
\section{Some relevant definitions}

It will be useful to recall two provably distinct combinatorial notions:   
slicewise VC dimension, and Shelah's $k$-dependence. For instance, when 
$\mcf$ is a family of bipartite graphs (subsets of $X \times Y$),  
\begin{itemize}
\item fixing any $a \in X$ induces a family of subsets of $Y$, 
$\{ \{ b \in Y : (a,b) \in F \} : F \in \mcf \}$, 
and likewise fixing any $b \in Y$ induces a family of subsets of $X$. 
Finite slicewise VC dimension means: there is some finite $n$ 
uniformly bounding the VC dimensions of all these families.  

\item for $k=2$, $2$-dependent means: for some finite $n$, every $\{ x_1, \dots, x_n \} \subseteq X$, 
every $\{ y_1, \dots, y_n \} \subseteq Y$, not all subsets of 
$\{ x_1, \dots, x_n \} \times \{ y_1, \dots, y_n \}$ arise as intersections with elements of $\mcf$. 
\end{itemize} 
We will also refer to a crucial difference in two definitions of statistical learning, 
which can be seen in the information the learner receives. 
For the comparison, let me again use the very simple example of learning a family 
$\mcf$ of bipartite graphs on $X \times Y$. 
\begin{itemize}
\item Kobayashi-Kuriyama-Takeuchi [K][KT]: the learner receives finitely many points 
$x_1, \dots, x_n$ from $X$ and $y_1, \dots, y_n$ from $Y$, along with the complete information 
of the neighborhoods of each of these points in the entire graph (which may be an infinite amount of 
information).  

\item Coregliano-Malliaris ``high-arity PAC learning'' [CM24]: the learner receives finitely 
many points $x_1, \dots, x_n$ from $X$ and $y_1, \dots, y_n$ from $Y$ along with the 
information of the induced subgraph \emph{just between these finitely many points}.
\end{itemize}
Coregliano-Malliaris proved that in both the partite and nonpartite case,  
 high-arity PAC learning characterizes finite slicewise VC dimension [CM24], [CM25a], and  
$k$-dependence is characterized by sample completion high-arity PAC learning, 
which I won't define here  
(a quite different learning model, about reconstruction on high-dimensional finite grids after random erasure) 
[CM25b].  

\section{Results of Chernikov-Towsner version 1}

In version 1 of Chernikov-Towsner (October 2, 2025), 
Theorem 6.4 says that ``every class $\mcf$ of finite $VC_k$-dimension is
$PAC_k$-learnable.''  
Per the definitions of dimension and learning in version 1, the assertion is    
that ``every $k$-dependent class is high-arity PAC learnable in the sense of [CM24].''

This is clearly false. A simple counterexample is the following. 
Suppose $X = Y = \mathbb{N}$ and work over $X \times Y$.  
Let $\mch$ be a family of subsets of $X$ which is not learnable in the classical PAC theory, 
for instance, $\mcp(X)$. Let $\mcf = \{ A \times Y : A \in \mch \}$. Clearly 
$\mcf$ is 2-dependent, as already we cannot shatter a $2\times 2$ square. On the other hand, 
it is not hard to see that existence of a high-arity PAC learner for $\mcf$ 
would yield a classic PAC learner for $\mch$, contradicting classical PAC theory. 
In order to avoid counterexamples of this kind, what we need is not 2-dependence, 
but rather finite slicewise VC dimension, which prevents a lifting of 
complexity from the lower dimensional case (as in [CM24]).    

One aside: in defining ``$PAC_k$-learnable'' in version 1,  
Chernikov-Towsner cite [K][KT], but give the definition from [CM24], as can clearly be seen      
from their page 5. They are aware there is a discrepancy, but they claim 
that the two definitions are equivalent using their Corollary 6.5, which 
depends on their false Theorem 6.4.

\section{Results of Chernikov-Towsner version 2} 

In version 2 of Chernikov-Towsner (October 15, 2025), what is now 
 Theorem 6.5 says that ``every class $\mcf$ of finite $VC_k$-dimension 
is properly $PAC_k$ learnable.'' 

``Proper'' means that the algorithm outputs an answer in the original class $\mcf$.
More importantly, the definition of learning has changed. They are now using the [K][KT] definition of learning, 
allowing unbounded slices. Per the definitions of dimension and learning in version 2, 
their theorem asserts that ``every $k$-dependent class is learnable in the sense of Kobayashi-Kuriyama-Takeuchi.'' 

At this point, of course, they are no longer proving results related to our work because 
they are using a completely different notion of learning. Specifically, it is unrelated to 
both high-arity PAC of [CM24] (which characterized finite slicewise VC dimension) 
and sample-completion high-arity PAC of [CM25b] (which characterized $k$-dependence 
by means of a bounded learning problem).  I also note that their packing statement for 
$k$-dependence in their Corollary 1.1  has no connection to either of the 
packing statements we had proved for the $k$-dependent case in [CM25b], which is 
to be expected since the learning is different.  

I have not looked into the correctness of the statement of their (v2) Theorem 6.5.  
However, the proof contains elementary mistakes. 
Their proposed proof is roughly structured as follows. First they use 
$k$-dependence to obtain a packing lemma.  
Now given a measure, 
the packing lemma gives centers $S_1, \dots, S_M$ with $M = M(\epsilon)$.  
By weak law of large numbers, if the sample size is large enough, they
claim that with sufficiently high probability,
 all empirical distances [i.e., distances on the sample] from the
unknown $S$ to each of the $S_i$ are close to the actual 
$\mu(S \Delta S_i)$. 
Their learning algorithm chooses the lexicographically first 
among the allowed Boolean combinations of the $S_i$'s which minimizes 
the empirical distance to $S$, and then does some additional work 
  to ensure properness (as their packing statement involves Boolean combinations and projections).  

Here is why this argument is incorrect: even though the packing bound $M(\epsilon)$ does not depend on the
underlying measure, the actual centers do.
Since the algorithm does not know the underlying measure, it cannot necessarily know the $S_i$.

\br
Finally, observe that Chernikov-Towsner make an unusual sequence of assumptions.  
They restrict only to finite spaces. Since any finite set is learnable given a large
enough sample, quantitative concerns would typically come into focus. 
However, they work in a completely qualitative setting, 
and use a definition of learning in which the learner receives a possibly unbounded 
amount of information.  Any one of these choices might 
make sense in isolation, but together they are puzzling.

\section{Minor points of contact}

With the errors and implicit retractions above, and since their notion of learning 
has changed, the significant interactions of [CT25]  
with work of Coregliano-Malliaris evaporate. Two minor points of contact are worth mentioning. 

First, [CT25] Remark 3.3 complains about our Lemma 7.8, stated across [CM25b] pages 45-46, where  
we derive and then apply a precise form of a combinatorial bound whose shape was first found by Shelah,
as we cite. We are glad to hear that a similar calculation was already done by 
Chernikov, Palacin and Takeuchi in the Notre Dame Journal [CPT] Proposition 3.9.  
Certainly, we are happy to add a citation.  
This was neither difficult combinatorics nor an intentional slight.  
I emphasize that the credit for $k$-dependence goes to Shelah. 

Second, Chernikov and Towsner claim that the main result of [CM25a] is a 
slicewise packing lemma and that this particular  
packing lemma is implicit in their earlier work [CT20]. 
They provide a somewhat vague proof relying heavily on [CT20]. 
Here are the parts of the statement I can respond to.  
(a) Until a readable, self-contained proof is provided, it is difficult to determine exactly what 
follows from [CT20]. 
(b) Because they work in product spaces, they only reproduce the ``partite'' version 
of this particular packing lemma, not the ``nonpartite'' one which was 
crucial for [CM25a] (think: independence property versus random graph).  
(c) The actual main results of [CM25a] per \S 5, Main Results, were to prove the implication that 
nonpartite high-arity PAC learning implies finite slicewise VC dimension, and the role of this packing lemma 
was to interpolate between these conditions, being as weak as possible modulo implying finite slicewise VC dimension
so as to be provable from nonpartite high-arity PAC.  Simply re-deriving that this weak packing statement is equivalent to finite 
slicewise VC dimension is not the main (or even a major) point.


\begin{thebibliography}{00}

\bibitem[CT25]{CT25} Artem Chernikov and Henry Towsner. ``Higher Arity PAC Learning, VC Dimension and Packing Lemma.'' 
ArXiv:2510.02420v1, October 2, 2025 and v2, October 15, 2025. 

\bibitem[CT]{CT} Artem Chernikov and Henry Towsner. ``Hypergraph regularity and higher arity
VC-dimension.'' Preprint, arXiv:2010.00726, 2020.

\bibitem[CPT]{CPT} Artem Chernikov, Daniel Palacin, and Kota Takeuchi. ``On n-dependence.''
Notre Dame Journal of Formal Logic, 60(2):195--214, 2019. (arXiv:1411.0120)

\bibitem[CM24]{CM24} Leonardo N. Coregliano and Maryanthe Malliaris, ``High-arity PAC learning via exchangeability.'' arXiv:2402.14294. 
145 pages. 

\bibitem[CM25a]{CM25a} Leonardo N. Coregliano and Maryanthe Malliaris, ``A packing lemma for $VCN_k$-dimension and learning
high-dimensional data.” arXiv:2505.15688. 29 pages.  

\bibitem[CM25b]{CM25b} Leonardo N. Coregliano and Maryanthe Malliaris, ``Sample completion, structured correlation 
and Netflix problems.'' arXiv:2509.20404. 97 pages. 

\bibitem[K]{K} Munehiro Kobayashi. ``A generalization of the PAC learning in product probability spaces.'' 
RIMS preprint, 2015.  
Available at https://www.kurims.kyoto-u.ac.jp/$\sim$kyodo/kokyoroku/contents/pdf/1938-06.pdf

\bibitem[KT]{KT} Takayuki Kuriyama and Kota Takeuchi. ``On the $PAC_n$ learning.'' RIMS preprint, 2015. 
Available at https://www.kurims.kyoto-u.ac.jp/$\sim$kyodo/kokyoroku/contents/pdf/1938-09.pdf

\bibitem[S]{S} Saharon Shelah. ``Strongly dependent theories.'' Israel J. Math. 204, 1 (2014),
pp. 1–83. 
\end{thebibliography}
\end{document}